\documentclass[11pt]{article}
\usepackage[textwidth=6.2in,textheight=8.5in,centering]{geometry}
\usepackage{microtype}
\usepackage{amsmath,amssymb,amsthm}
\usepackage{graphicx,tikz}
\usetikzlibrary{arrows.meta}
\usepackage{enumitem,float,authblk}
\usepackage{cite}
\usepackage[hidelinks]{hyperref}

\newtheorem{theorem}{Theorem}[section]
\newtheorem{lemma}[theorem]{Lemma}
\newtheorem{claim}{Claim}

\newcommand{\chic}{\chi_{\mathrm c}}
\newcommand{\chib}{\chi_{\mathrm b}}
\newcommand{\R}{\mathbb R}
\newcommand{\Z}{\mathbb Z}
\title{A planar signed graph with circular chromatic number five}
\author[1]{Lujia Wang} \author[2]{Huan Zhou}

\affil[1]{\small Hainan Bielefeld University of Applied Sciences, Danzhou, China. {Email: \texttt{lujia.wang@hibiuh.edu.cn}}}
\affil[2]{\small School of Science, Jimei University, Xiamen, China. {Email: \texttt{huanzhou@zjnu.edu.cn}}} 

\date{}
\begin{document} \maketitle 

\begin{abstract}
A circular $r$-coloring of a signed graph maps its vertices to a circle
of circumference $r\ge2$, with the endpoints of negative edges at circular
distance at least $1$ and those of positive edges at most $r/2-1$;
the circular chromatic number is the infimum of such $r$.
We construct a simple planar signed triangulation on $85$ vertices with
circular chromatic number exactly $5$. This improves the lower bound on
the supremum of the circular chromatic numbers of simple planar signed
graphs from $14/3$ to $5$, giving a partial answer to a question posed
by Naserasr, Wang, and Zhu.
\end{abstract}

\section{Introduction}\label{sec:introduction}

A signed graph $(G,\sigma)$ consists of a graph $G$ and a signature
$\sigma:E(G)\to\{+1,-1\}$. Zaslavsky's zero-free $k$-coloring assigns
colors from $\{\pm1,\ldots,\pm k\}$ so that
$f(u)\ne\sigma(uv)f(v)$ on every edge~\cite{Z1982Coloring}.
Reversing all edge signs gives an equivalent formulation in terms of
balanced vertex sets. A set is \emph{balanced} if every cycle in its
induced signed subgraph has positive edge-sign product. The least
number of balanced sets partitioning $V(G)$ is the \emph{balanced
chromatic number} $\chib(G,\sigma)$, called the \emph{balanced partition
number} by Zaslavsky~\cite{Z1987}. Thus a balanced $k$-coloring of
$(G,\sigma)$ is equivalent to a zero-free coloring of $(G,-\sigma)$
using $k$ signed pairs.

Jim\'enez, McDonald, Naserasr, Nurse, and Quiroz~\cite{JMNNQ2026}
reintroduced this viewpoint as balanced coloring in their study of
Hadwiger-type questions. Related work on fractional
balanced coloring~\cite{NPPZ2026} and signed analogues of Kneser,
Schrijver, and Borsuk graphs~\cite{KNWYZZ2025} connects this parameter
with fractional and topological methods. Here we study its circular
refinement, with particular attention to planar graphs.

For unsigned graphs, circular coloring places vertices on a circle
with adjacent vertices at distance at least one. The least possible
circumference is the circular chromatic number. Introduced by
Vince~\cite{V1988} as the star chromatic number, it refines ordinary
coloring by allowing nonintegral values; see Zhu's
survey~\cite{Zhu2001}. The signed version studied here is due to
Naserasr, Wang, and Zhu~\cite{NWZ2021}.

For $r\ge2$, a \emph{circular $r$-coloring} of $(G,\sigma)$ is a map
$\phi:V(G)\to\R/r\Z$ such that
\[
\begin{aligned}
 d_r(\phi(u),\phi(v))&\ge1 &&\text{if }\sigma(uv)=-1,\\
 d_r(\phi(u),\phi(v))&\le r/2-1 &&\text{if }\sigma(uv)=+1,
\end{aligned}
\]
where $d_r$ denotes distance along the shorter arc. The \emph{circular
chromatic number} $\chic(G,\sigma)$ is the infimum of such $r$.
A balanced $k$-coloring exists if and only if a circular $2k$-coloring
exists~\cite{NPW2026}; thus circular coloring refines balanced coloring.

This definition reverses the edge signs in the convention of
Naserasr, Wang, and Zhu~\cite{NWZ2021}; it agrees with the convention
used by Chen, Naserasr, and Sarti~\cite{CNS2026}. We state all results
in this notation. With all edges negative, the definition reduces to
ordinary circular coloring.

Recent work has developed both structural and topological approaches
to this parameter. Naserasr and Zhou~\cite{NZ2025} proved the sharp
bound $10/3$ for connected simple signed graphs of maximum degree at
most $3$, apart from the switching class of the all-negative $K_4$.
Li, Naserasr, Wang, and Zhu~\cite{LNWZ2024} developed the corresponding
circular-flow theory and its duality with coloring for plane signed
graphs. Naserasr, Pujol, and
Wang~\cite{NPW2026} developed a unified winding-number framework for
circular-coloring lower bounds and gap theorems, including results for
unsigned and signed graphs embedded on surfaces.
%Winding-number arguments were used by Gujgiczer, Naserasr, Rohini, and Taruni~\cite{GNST2023} to construct signed bipartite graphs with circular chromatic number $4$.

Let $\mathcal{SP}$ denote the class of simple planar signed graphs, and
write
\[
 \chic(\mathcal{SP})=
 \sup\{\chic(G,\sigma):(G,\sigma)\in\mathcal{SP}\}.
\]
Kardo\v{s} and Narboni~\cite{KN2021} constructed simple planar signed
graphs with circular chromatic number greater than $4$.
Naserasr, Wang, and Zhu~\cite{NWZ2021} subsequently established
\[
 \frac{14}{3}\le\chic(\mathcal{SP})\le6
\]
and asked for the exact supremum. Their lower bound comes from an
explicit example with circular chromatic number $14/3$; the upper
bound follows by greedy coloring with six integer colors modulo $6$:
each colored neighbor forbids just one color, and every simple planar
graph is $5$-degenerate. Gu and Zhu~\cite{GZ2022} further showed that every rational
value in $[2,14/3]$ occurs in $\mathcal{SP}$.

Sparsity gives substantially smaller bounds. Simple signed
$2$-degenerate graphs and simple signed bipartite planar graphs have
circular chromatic number strictly less than $4$~\cite{KNNW2023}.
For signed planar graphs, the same holds when the underlying graph
has girth at least $5$~\cite{LLW2025}; girth at least $6$ gives a bound
of $3$~\cite{BHNSW2024}, and girth at least $8$ gives a strict bound
below $3$~\cite{LMSWW2026}. These results contrast with the larger
values possible in the unrestricted planar class.

For unrestricted simple planar signed graphs, we improve the lower
bound from $14/3$ to $5$.
\begin{theorem}\label{thm:circular}
There is a simple plane signed triangulation $(G,\sigma)$ on $85$
vertices with $\chic(G,\sigma)=5$.
\end{theorem}

In Section~\ref{sec:construction}, we give the construction and prove
Theorem~\ref{thm:circular}. Section~\ref{sec:concluding} presents the
frame of a smaller example and discusses the remaining open problem.

\section{The construction and its circular chromatic number}
\label{sec:construction}

\subsection{Two triangular fillings}\label{subsec:fillings}

A cycle of a signed graph is \emph{positive} or \emph{negative}
according to the product of its edge signs. \emph{Switching} a
vertex reverses the signs of all its incident edges. Every cycle
uses zero or two such edges, so switching preserves cycle signs.
Two signed graphs on the same underlying graph are \emph{equivalent}
if one can be obtained from the other by a sequence of switchings.
By Zaslavsky's characterization~\cite[Proposition~3.2]{Z1982}, this
holds if and only if every cycle has the same sign in both graphs.
In a circular $r$-coloring, switching $v$ is realized by replacing
$\phi(v)$ with $\phi(v)+r/2$, since
\[
 d_r(\phi(u),\phi(v)+r/2)=r/2-d_r(\phi(u),\phi(v)).
\]
Thus equivalent signed graphs have the same circular chromatic number.

Throughout the lower-bound argument, fix $4\le r<5$ and a circular
$r$-coloring $\phi$, and abbreviate $d_r(\phi(u),\phi(v))$ to $d(u,v)$.
The restriction $r\ge4$ will be justified by the all-negative $K_4$
contained in our construction.
In a fixed circular $r$-coloring, a negative edge $uv$ is \emph{short} if $d(u,v)<2$, and a positive
edge $uv$ is \emph{short} if $d(u,v)>r/2-2$. An edge that is not
short is \emph{long}. By the switching identity, either condition
says that the endpoint colors have distance less than $2$ when the
edge is made negative. Thus switching preserves short and long edges.

First observe that a negative triangle has at most one long edge.
Indeed, switch it to all-negative signs. Its three colors leave three
circular gaps, each of length at least $1$. Two long edges would give
two gaps of length at least $2$, forcing $r\ge2+2+1=5$.

Filling a triangular face can further restrict the long edges on its
boundary. Figure~\ref{fig:small-fillings} shows two signed plane graphs
$T^-$ and $T^+$, with negative and positive outer triangles, respectively.
To fill a face, identify the outer triangle of the chosen graph with
the face boundary and place its remaining vertices inside the face.
Switch the copy to match the boundary signs; this is possible whenever
the two triangles have the same sign.

\begin{figure}[H]
\centering
\begin{minipage}[c]{.48\textwidth}
\centering
\begin{tikzpicture}[x=2.1cm,y=2.1cm,
  every node/.style={font=\footnotesize,inner sep=.8pt},
  negative/.style={draw=red,line width=.8pt},
  vertex/.style={circle,draw=black!65,fill=white,line width=.35pt,minimum size=12pt}]
\coordinate (A) at (0,1.5);
\coordinate (B) at (-1.3,-.8);
\coordinate (C) at (1.3,-.8);
\coordinate (t) at (0,-.033333);
\draw[negative] (A)--(B)--(C)--cycle;
\draw[negative] (A)--(t)--(B);
\draw[negative] (C)--(t);
\node[vertex] at (A) {$u$};
\node[vertex] at (B) {$v$};
\node[vertex] at (C) {$w$};
\node[vertex] at (t) {$t$};
\end{tikzpicture}
\par\smallskip\footnotesize (a) $T^-$
\end{minipage}
\hfill
\begin{minipage}[c]{.48\textwidth}
\centering
\begin{tikzpicture}[x=2.1cm,y=2.1cm,
  every node/.style={font=\footnotesize,inner sep=.8pt},
  positive/.style={draw=blue,dash pattern=on 3pt off 2pt,line width=.8pt},
  negative/.style={draw=red,line width=.8pt},
  vertex/.style={circle,draw=black!65,fill=white,line width=.35pt,minimum size=12pt}]
\coordinate (va) at (-1.300000,-0.800000);
\coordinate (vb) at (0.000000,1.500000);
\coordinate (vc) at (1.300000,-0.800000);
\coordinate (vx) at (0.000000,-0.368196);
\coordinate (vy) at (-0.290000,0.134098);
\coordinate (vz) at (0.290000,0.134098);
\coordinate (vda) at (-0.530000,-0.344699);
\coordinate (vdb) at (0.000000,0.589399);
\coordinate (vdc) at (0.530000,-0.344699);
\draw[positive] (va)--(vb)--(vc)--cycle;
\draw[positive] (va)--(vy);
\draw[positive] (vb)--(vz);
\draw[positive] (vc)--(vx);
\draw[negative] (va)--(vx);
\draw[negative] (vb)--(vy);
\draw[negative] (vc)--(vz);
\draw[negative] (vx)--(vy)--(vz)--cycle;
\draw[positive] (va)--(vda)--(vx);
\draw[negative] (vda)--(vy);
\draw[positive] (vb)--(vdb);
\draw[positive] (vdb)--(vy);
\draw[negative] (vdb)--(vz);
\draw[positive] (vc)--(vdc)--(vz);
\draw[negative] (vdc)--(vx);
\node[vertex] at (va) {$v$};
\node[vertex] at (vb) {$u$};
\node[vertex] at (vc) {$w$};
\node[vertex] at (vx) {$x$};
\node[vertex] at (vy) {$y$};
\node[vertex] at (vz) {$z$};
\node[vertex] at (vda) {$t_v$};
\node[vertex] at (vdb) {$t_u$};
\node[vertex] at (vdc) {$t_w$};
\end{tikzpicture}
\par\smallskip\footnotesize (b) $T^+$
\end{minipage}

\caption{The signed plane graphs $T^-$ and $T^+$, with all-negative and
all-positive outer triangles $uvw$, respectively.
Positive edges are blue and dashed; negative edges are red and solid.} 
\label{fig:small-fillings}
\end{figure}
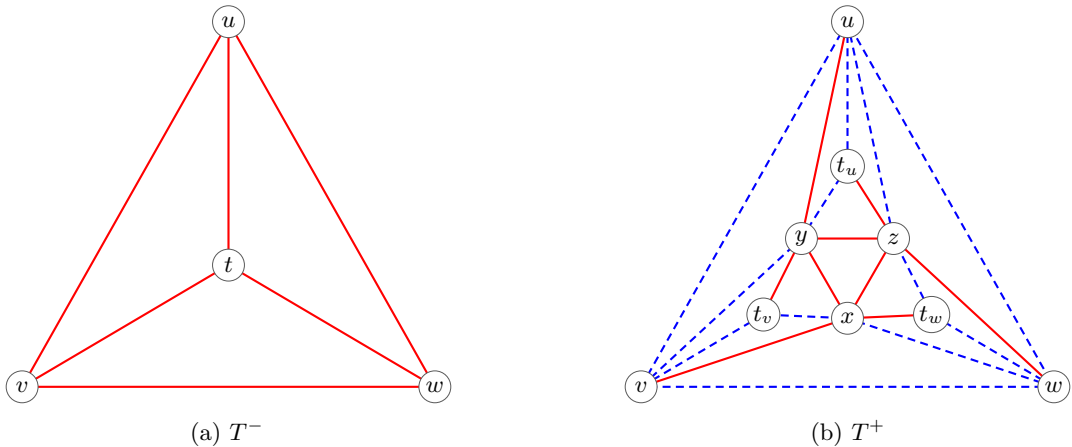
\begin{lemma}\label{lem:filling}
A signed triangular face admits a filling such that, in every circular
$r$-coloring of the resulting graph with $4\le r<5$:
\begin{enumerate}[label=\textup{(\roman*)}]
\item if its boundary is negative, exactly one boundary edge is long;
\item if its boundary is positive, at most one boundary edge is long.
\end{enumerate}
The graphs $T^-$ and $T^+$ in Figure~\ref{fig:small-fillings} give the
respective fillings.
\end{lemma}
\begin{proof}
For \textup{(i)}, use $T^{-}=(K_4,-)$. Its four colors leave
four gaps, each of length at least $1$. Deleting the interior vertex's
color merges two consecutive gaps. Both the merged arc and its
complement have length at least $2$, while each of the two remaining
gaps has length at most $r-3<2$. Thus exactly one boundary edge is long.

For \textup{(ii)}, put $\delta=r/2-2$ and $\ell=r/2-1=1+\delta$, so
$0\le\delta<1/2$. A positive edge has distance at most $\ell$, and it is
long when its distance is at most $\delta$. By cyclic symmetry, it
suffices to rule out two long edges $uv,uw$.
Set $\phi(u)=0$, use representatives in $[-r/2,r/2)$, and abbreviate
colors by vertex names. Then $v,w\in[-\delta,\delta]$.
The negative edge $vx$ excludes $x=0$; reflecting the circle if necessary,
assume $x<0$. The edges $wx,vx$ give
\[
 w-\ell\le x\le v-1.
\]
If $y<0$, the edges $uy,vy$ give $y\in[v-\ell,-1]$; if $z<0$,
the edges $uz,wz$ give $z\in[-\ell,w-1]$. Either interval puts that color within ordinary distance
at most $2\delta<1$ of $x$, contradicting $xy$ or $xz$.
Thus $y,z>0$, and
\[
 1\le y\le v+\ell<2,\qquad w+1\le z\le\ell.
\]
Since $z-y\le\ell-1=\delta<1$, the negative edge $yz$ forces $y\ge z+1$.

Write $t$ for the color of $t_u$. The positive edges $ut_u,yt_u$
restrict $t$ to $[y-\ell,\ell]$. There is no overlap around the other
side of the circle, since $y<2=r-2\ell$.
But $y-z\ge1$ and $z\ge w+1\ge1-\delta$, so
\[
 [y-\ell,\ell]\subseteq[z-\delta,z+2\delta]
 \subset(z-1,z+1).
\]
The last inclusion is strict because $0\le\delta<1/2$.
Thus every available color for $t_u$ is at distance less than $1$
from $z$, contradicting the negative edge $zt_u$.
Only $t_u$ is needed for this pair of boundary edges; the other two
added vertices play the same role for the other two pairs.
\end{proof}

\subsection{The construction}\label{subsec:construction}

We first construct a signed plane graph $\widehat H$ from the base in
Figure~\ref{fig:H-construction}, with outer triangle $w_1w_2w_3$.
Fill each red triangle with a copy of $T^-$
and each blue triangle with a copy of $T^+$. In each copy of $T^+$,
switch $u$ and then identify $w$ with $s$ and $u$ with the corner
opposite the positive edge. The boundary signs then agree with those
of the blue triangle. Thus $|V(\widehat H)|=11+7+2\cdot6=30$.

\begin{figure}[htbp]
\centering
\begin{minipage}[c]{.53\textwidth}\centering
\begin{tikzpicture}[x=.65cm,y=.65cm,
  every node/.style={font=\scriptsize,inner sep=1pt},
  positive/.style={draw=blue,dash pattern=on 3pt off 2pt,line width=.8pt},
  negative/.style={draw=red,line width=.8pt},
  vertex/.style={circle,draw=black!65,fill=white,line width=.35pt,minimum size=11pt}]
\coordinate (a) at (-4.734272,0.000000);
\coordinate (b) at (0.000000,8.200000);
\coordinate (c) at (4.734272,0.000000);
\coordinate (d) at (-1.050000,4.551987);
\coordinate (e) at (1.050000,4.551987);
\coordinate (f) at (-2.100000,2.733333);
\coordinate (g) at (1.050000,0.914680);
\coordinate (h) at (2.100000,2.733333);
\coordinate (i) at (-1.050000,0.914680);
\coordinate (j) at (-0.350000,3.339551);
\coordinate (k) at (0.350000,2.127116);
\fill[red!4] (a)--(b)--(d)--cycle;
\fill[red!4] (b)--(c)--(h)--cycle;
\fill[red!4] (b)--(e)--(h)--cycle;
\fill[blue!4] (d)--(j)--(f)--cycle;
\fill[blue!4] (e)--(k)--(j)--cycle;
\fill[red!4] (e)--(h)--(k)--cycle;
\fill[red!4] (g)--(h)--(k)--cycle;
\fill[red!4] (c)--(g)--(h)--cycle;
\fill[red!4] (a)--(i)--(g)--cycle;
\draw[negative] (a)--(b);
\draw[negative,line width=1.15pt] (a)--(c);
\draw[negative] (a)--(d);
\draw[negative] (a)--(f);
\draw[negative] (a)--(g);
\draw[negative] (a)--(i);
\draw[negative] (b)--(c);
\draw[negative] (b)--(d);
\draw[negative] (b)--(e);
\draw[negative] (b)--(h);
\draw[negative] (c)--(g);
\draw[negative] (c)--(h);
\draw[negative] (d)--(e);
\draw[negative] (d)--(f);
\draw[positive] (d)--(j);
\draw[negative] (e)--(h);
\draw[positive] (e)--(j);
\draw[negative] (e)--(k);
\draw[negative] (f)--(i);
\draw[negative] (f)--(j);
\draw[negative] (g)--(h);
\draw[negative] (g)--(i);
\draw[negative] (g)--(k);
\draw[negative] (h)--(k);
\draw[negative] (i)--(j);
\draw[negative] (i)--(k);
\draw[negative] (j)--(k);
\node[vertex] at (a) {$w_1$};
\node[vertex] at (b) {$w_3$};
\node[vertex] at (c) {$w_2$};
\node[vertex] at (d) {$s_1$};
\node[vertex] at (e) {$s_2$};
\node[vertex] at (f) {$t$};
\node[vertex] at (g) {$t_3$};
\node[vertex] at (h) {$z$};
\node[vertex] at (i) {$t_1$};
\node[vertex] at (j) {$s$};
\node[vertex] at (k) {$t_2$};
\node[font=\tiny,text=black!65] at (-3.30000,5.80000) {$O$};
\node[font=\tiny,text=black!65] at (0.35000,0.30489) {$V$};
\node[font=\tiny,text=black!65] at (-1.92809,4.25066) {$R$};
\node[font=\tiny,text=black!65] at (-1.57809,0.60979) {$L$};
\node[font=\tiny,text=black!65] at (2.27809,3.64444) {$E$};
\node[font=\tiny,text=black!65] at (2.62809,1.21600) {$F$};
\node[font=\tiny,text=black!65] at (1.05000,5.16177) {$H$};
\node[font=\tiny,text=black!65] at (1.16667,1.92504) {$K$};
\node[font=\tiny,text=black!65] at (0.00000,5.76799) {$X$};
\node[font=\tiny,text=black!65] at (0.11667,1.31883) {$Y$};
\node[font=\tiny,text=black!65] at (-0.11667,4.14784) {$A$};
\node[font=\tiny,text=black!65] at (-0.35000,2.12712) {$B$};
\node[font=\tiny,text=black!65] at (-2.62809,2.42844) {$C$};
\node[font=\tiny,text=black!65] at (-2.62809,1.21600) {$D$};
\node[font=\tiny,text=black!65] at (-1.16667,2.32919) {$M$};
\node[font=\tiny,text=black!65] at (1.16667,3.13748) {$U$};
\node[font=\tiny,text=black!65] at (-1.16667,3.54162) {$P$};
\node[font=\tiny,text=black!65] at (0.35000,3.33955) {$Q$};
\end{tikzpicture}
\par\smallskip\footnotesize (a) The signed base of $\widehat H$.
\end{minipage}\hfill
\begin{minipage}[c]{.44\textwidth}\centering
\begin{tikzpicture}[x=.65cm,y=.65cm,
 every node/.style={font=\scriptsize,inner sep=1pt},
 bare/.style={circle,draw=black!60,fill=white,minimum size=14pt,line width=.4pt},
 filled/.style={bare,draw=red!70!black,fill=red!8},
 plus/.style={rectangle,rounded corners=2pt,draw=blue!70!black,fill=blue!8,minimum size=14pt}]
\coordinate (rR) at (-3,2.5);
\coordinate (rX) at (-1,2.5);
\coordinate (rH) at (1,2.5);
\coordinate (rE) at (2.5,2.5);
\coordinate (rL) at (-3,-2.5);
\coordinate (rY) at (-1,-2.5);
\coordinate (rK) at (1,-2.5);
\coordinate (rF) at (2.5,-2.5);
\coordinate (rO) at (-3,4);
\coordinate (rV) at (-3,-4);
\coordinate (rC) at (-3,1);
\coordinate (rD) at (-3,-1);
\coordinate (rP) at (-2,1);
\coordinate (rM) at (-2,-1);
\coordinate (rA) at (-1,1);
\coordinate (rB) at (-1,-1);
\coordinate (rQ) at (0,0);
\coordinate (rU) at (1,0);
\draw[black,line width=1.25pt] (rO)--(-4.1,4)--(-4.1,-4)--(rV);
\draw[black!45,line width=.5pt] (rO)--(rR);
\draw[black!45,line width=.5pt] (rO)--(2.5,4)--(rE);
\draw[black!45,line width=.5pt] (rV)--(rL);
\draw[black!45,line width=.5pt] (rV)--(2.5,-4)--(rF);
\draw[black!45,line width=.5pt] (rR)--(rX);
\draw[black!45,line width=.5pt] (rR)--(rC);
\draw[black!45,line width=.5pt] (rL)--(rY);
\draw[black!45,line width=.5pt] (rL)--(rD);
\draw[black!45,line width=.5pt] (rE)--(rF);
\draw[black!45,line width=.5pt] (rE)--(rH);
\draw[black!45,line width=.5pt] (rF)--(rK);
\draw[black!45,line width=.5pt] (rH)--(rX);
\draw[black!45,line width=.5pt] (rH)--(rU);
\draw[black!45,line width=.5pt] (rK)--(rY);
\draw[black!45,line width=.5pt] (rK)--(rU);
\draw[black!45,line width=.5pt] (rX)--(rA);
\draw[black!45,line width=.5pt] (rY)--(rB);
\draw[black!45,line width=.5pt] (rA)--(rP);
\draw[black!45,line width=.5pt] (rA)--(rQ);
\draw[black!45,line width=.5pt] (rB)--(rM);
\draw[black!45,line width=.5pt] (rB)--(rQ);
\draw[black!45,line width=.5pt] (rC)--(rD);
\draw[black!45,line width=.5pt] (rC)--(rP);
\draw[black!45,line width=.5pt] (rD)--(rM);
\draw[black!45,line width=.5pt] (rM)--(rP);
\draw[black!45,line width=.5pt] (rU)--(rQ);
\node[filled] at (rR) {$R$};
\node[bare] at (rX) {$X$};
\node[filled] at (rH) {$H$};
\node[filled] at (rE) {$E$};
\node[filled] at (rL) {$L$};
\node[bare] at (rY) {$Y$};
\node[filled] at (rK) {$K$};
\node[filled] at (rF) {$F$};
\node[bare] at (rO) {$O$};
\node[bare] at (rV) {$V$};
\node[bare] at (rC) {$C$};
\node[bare] at (rD) {$D$};
\node[plus] at (rP) {$P$};
\node[bare] at (rM) {$M$};
\node[bare] at (rA) {$A$};
\node[bare] at (rB) {$B$};
\node[plus] at (rQ) {$Q$};
\node[filled] at (rU) {$U$};
\end{tikzpicture}
\par\smallskip\footnotesize (b) The same regions, arranged for the proof.
\end{minipage}
\caption{The $30$-vertex gadget $\widehat H$. In (a), positive edges
are blue dashed and negative edges are red solid. Fill red regions with
$T^-$ and blue regions $P,Q$ with $T^+$. In (b), nodes represent regions:
red circles indicate $T^-$ fillings, blue boxes indicate $T^+$ fillings,
and white circles indicate unfilled negative triangles. The unbounded
region is $O$. Lines record shared boundary edges; the thick line $OV$
represents $w_1w_2$.}
\label{fig:H-construction}
\end{figure}

\begin{figure}[H]
\centering
\begin{tikzpicture}[x=2cm,y=2cm,
  every node/.style={font=\footnotesize,inner sep=1pt},
  negative/.style={draw=red,line width=.8pt},
  vertex/.style={circle,draw=black!65,fill=white,line width=.35pt,minimum size=13pt}]
\coordinate (a) at (0,1.5);
\coordinate (b) at (-1.3,-.8);
\coordinate (c) at (1.3,-.8);
\coordinate (d) at (0,-.03);
\draw[negative,line width=1.2pt] (a)--(b)--(c)--cycle;
\draw[negative] (a)--(d)--(b);
\draw[negative] (c)--(d);
\node at (-.43,.22) {$\widehat{H}_1$};
\node at (0,-.54) {$\widehat{H}_2$};
\node at (.43,.22) {$\widehat{H}_3$};
\node[vertex] at (a) {$u$};
\node[vertex] at (b) {$v$};
\node[vertex] at (c) {$x$};
\node[vertex] at (d) {$y$};
\end{tikzpicture}
\caption{The construction of $(G,\sigma)$ from $(K_4,-)$.
The thick edges $uv,vx,xu$ correspond to $w_1w_2$ in the three copies
of $\widehat H$.}
\label{fig:primal-construction}
\end{figure}
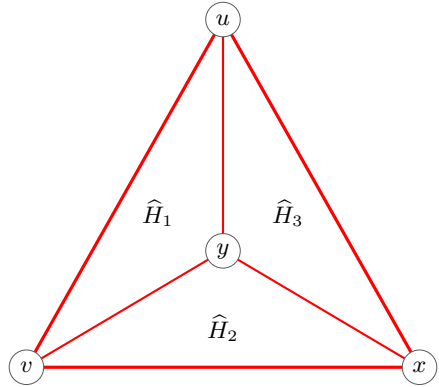

Now take an all-negative $K_4$ on $u,v,x,y$ with outer triangle $uvx$.
As shown in Figure~\ref{fig:primal-construction}, insert
$\widehat H_1,\widehat H_2,\widehat H_3$ into $uvy,vxy,xuy$, respectively,
identifying their boundary vertices by
\[
 (w_1,w_2,w_3)=(u,v,y),\quad(v,x,y),\quad(x,u,y),
\]
respectively. The resulting $(G,\sigma)$ is a simple plane triangulation with
\[
 |V(G)|=4+3\cdot(30-3)=85,\qquad |E(G)|=3\cdot85-6=249.
\]

By Lemma~\ref{lem:filling}(i), any circular $r$-coloring with $4\le r<5$
of its $K_4$ would give one of $uv,vx,xu$ distance at least $2$.
Each of these edges corresponds to $w_1w_2$ in a copy of $\widehat H$.
Thus the lower bound reduces to the following property of $\widehat H$.

\subsection{The marked edge}\label{subsec:marked-edge}

\begin{lemma}\label{circ:S}
In every circular $r$-coloring of $\widehat H$ with $4\le r<5$,
the edge $w_1w_2$ is short.
\end{lemma}

To prove the lemma, we compare the circular orders of the negative
triangular regions. First switch a region's boundary to all-negative,
and read its vertices counterclockwise around the region (clockwise
for the unbounded region). Their colors occur either clockwise or
counterclockwise on the color circle; this is the region's
\emph{circular order}. It is well defined, since the two switchings
making a triangle all-negative differ by shifting all three colors
by $r/2$.

\begin{claim}\label{clm:circular-orders}
In a circular $r$-coloring with $4\le r<5$, the following hold.
\begin{enumerate}[label=\textup{(\roman*)}]
\item Two negative triangular regions lying on opposite sides of a common short edge have opposite circular orders. Consequently, if their orders are the same, their common edge is long.
\item Suppose that each boundary edge of a copy of $T^+$ has a negative
triangular region on its other side. These three regions cannot all have the same circular order.
\item Under the hypothesis of \textup{(ii)}, if one boundary edge of $T^+$ is long, the negative triangular regions across its other two
boundary edges have opposite circular orders.
\end{enumerate}
\end{claim}

\begin{figure}[htbp]
  \centering
  \resizebox{\linewidth}{!}{%

\begin{tikzpicture}[x=1cm,y=1cm,line join=round,line cap=round,
  positive/.style={draw=blue,line cap=butt,dash pattern=on 3pt off 2pt,line width=.8pt},
  negative/.style={draw=red,line width=.8pt},
  order/.style={-{Stealth[length=2.4pt,width=2.1pt]},line width=.45pt,draw=black!85},
  vertex/.style={circle,draw=black!65,fill=white,line width=.35pt,
                 inner sep=0pt,minimum size=3.5pt},
  named vertex/.style={vertex,minimum size=12pt,inner sep=1pt},
  every node/.style={font=\small}]
\begin{scope}
  \coordinate (top) at (0,1.15);
  \coordinate (bottom) at (0,-1.65);
  \coordinate (left) at (-1.9,-.25);
  \coordinate (right) at (1.9,-.25);
  \draw[negative] (top)--(left)--(bottom)--(right)--cycle;
  \draw[negative] (top)--(bottom)
    node[midway,sloped,fill=white,inner sep=1pt,font=\scriptsize] {short};
  \node at (-.78,-.10) {$\Delta_1$};
  \node at (.78,-.10) {$\Delta_2$};
  \draw[order] (-.78,-.65) ++(55:.22) arc (55:355:.22);
  \draw[order] (.78,-.65) ++(125:.22) arc (125:-175:.22);
  \foreach \p in {top,bottom,left,right} \node[vertex] at (\p) {};
  \node at (0,-2.14) {(i)};
  \node at (0,-2.52) {Opposite orders};
\end{scope}
% (ii), (iii): T^+ has an all-positive boundary. Each exterior
% triangle has two positive edges and one negative edge.
\foreach \panel/\offset in {2/4.9,3/9.8}{
\begin{scope}[xshift=\offset cm]
  \coordinate (u) at (0,1.135641);
  \coordinate (v) at (-.8,-.25);
  \coordinate (w) at (.8,-.25);
  \coordinate (a) at (-1.6,1.135641);
  \coordinate (b) at (1.6,1.135641);
  \coordinate (c) at (0,-1.635641);
  \draw[negative] (u)--(a);
  \draw[negative] (w)--(b);
  \draw[negative] (v)--(c);
  \draw[positive] (v)--(a);
  \draw[positive] (u)--(b);
  \draw[positive] (w)--(c);
  \draw[positive] (v)--(u)
    node[midway,sloped,fill=white,inner sep=1pt,font=\scriptsize] {short};
  \draw[positive] (u)--(w)
    node[midway,sloped,fill=white,inner sep=1pt,font=\scriptsize] {short};
  \ifnum\panel=2
    \draw[positive] (v)--(w)
      node[midway,fill=white,inner sep=1pt,font=\scriptsize] {short};
  \else
    \draw[positive,line width=1.8pt] (v)--(w);
    \node[below=2pt,font=\footnotesize] at (0,-.25) {long};
  \fi
  \node at (0,.2) {$T^+$};
  \node at (-.8,.86) {$\Delta_1$};
  \node at (.8,.86) {$\Delta_2$};
  \draw[order] (-.8,.46) ++(55:.18) arc (55:355:.18);
  \draw[order] (.8,.46) ++(125:.18) arc (125:-175:.18);
  \ifnum\panel=2
    \node at (0,-.66) {$\Delta_3$};
    \draw[order] (0,-1.06) ++(55:.16) arc (55:355:.16);
  \else
    \node at (0,-1.03) {$\Delta_3$};
  \fi
  \ifnum\panel=2
    \node at (0,-2.14) {(ii)};
    \node at (0,-2.52) {Not all the same order};
  \else
    \node at (0,-2.14) {(iii)};
    \node at (0,-2.52) {$\Delta_1,\Delta_2$: opposite orders};
  \fi
  \foreach \p in {a,b,c} \node[vertex] at (\p) {};
  \foreach \p in {u,v,w} \node[named vertex] at (\p) {$\p$};
\end{scope}
}
\end{tikzpicture}
}
  \caption{Circular-order constraints. In \textup{(ii)} and \textup{(iii)},
  $u,v,w$ occur counterclockwise on the boundary of $T^+$, with
  $\phi(u)\le\phi(v)\le\phi(w)$ on their containing arc.
  Panel \textup{(ii)} shows the case in which all three boundary edges
  are short; panel \textup{(iii)} shows the case in which $vw$ is long.
  Arrows indicate the circular orders of the exterior negative triangles.}
  \label{fig:circular-order-rules}
\end{figure}
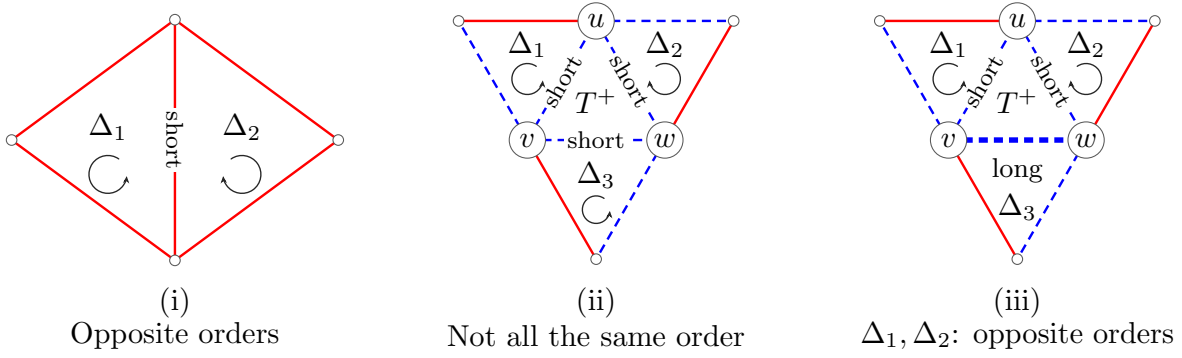

\begin{proof}
For \textup{(i)}, switch both triangles to all-negative. Neither third
color can lie on the shorter arc of their common short edge: this arc
has length less than $2$, while the third color must be at least $1$
from both endpoints. The two boundaries traverse the common edge in
opposite directions, so their circular orders are opposite.

For \textup{(ii)} and \textup{(iii)}, switch the $T^+$ boundary to
all-positive. Its three colors lie in an arc of length at most $r/2-1$.
Indeed, their largest gap has length at least $r/3>r/2-1$; since its
endpoints are joined by a positive edge, the complementary arc must
have length at most $r/2-1$.
Label the boundary $u,v,w$ counterclockwise. By reflecting the color
circle and cyclically relabeling if necessary, we may assume that
$\phi(u)\le\phi(v)\le\phi(w)$ along this arc in the counterclockwise
direction. The edge $uw$ has the greatest endpoint distance, so it
must be short: if it were long, all three boundary edges would be
long, contrary to Lemma~\ref{lem:filling}(ii). That lemma also ensures
that at least one of $uv,vw$ is short.

Let $\Delta_1,\Delta_2,\Delta_3$ be the exterior negative triangles
across $uv,uw,vw$, respectively, as shown in
Figure~\ref{fig:circular-order-rules}. Their orders follow from one
observation. Consider any short boundary edge with colors $a<b$ on
the containing arc. Switch the endpoint colored $a$ and, if necessary,
the third vertex of its exterior triangle to make that triangle
all-negative. The counterclockwise arc from $b$ to $a+r/2$ has length
\[
 r/2-(b-a)<2.
\]
As in \textup{(i)}, the third color lies outside this short arc.
Thus the exterior triangle has counterclockwise circular order if
its boundary traverses the shared edge from the later color to the
earlier one, and clockwise order if it traverses it in the reverse
direction. We apply this observation to each exterior triangle
separately.

The boundaries of $\Delta_1,\Delta_2,\Delta_3$ traverse their shared
edges as $vu,uw,wv$, respectively. Hence $\Delta_2$ has clockwise
order, while $\Delta_1$ and $\Delta_3$ have counterclockwise order
whenever $uv$ and $vw$, respectively, are short.
Since $uw$ and at least one of $uv,vw$ are short, two exterior
triangles have opposite orders, proving \textup{(ii)}.
If one boundary edge is long, it is $uv$ or $vw$; the other two are
short and their exterior triangles have opposite orders, proving
\textup{(iii)}.
\end{proof}

\begin{proof}[Proof of Lemma~\ref{circ:S}]
Use the region labels in Figure~\ref{fig:H-construction}(b).
For adjacent regions $S$ and $T$, $ST$ denotes their shared boundary edge;
$S=T$ and $S\ne T$ mean that their circular orders are equal and
opposite, respectively. First note the central obstruction:
\[A=C, \qquad B=D \qquad \implies \qquad A=B.\]
Indeed, Claim~\ref{clm:circular-orders}(ii) at $P$ forces $M$ to
have the opposite order to $A=C$. If $B=D\ne A$, then $B=D=M$,
so Claim~\ref{clm:circular-orders}(i) makes both $MB$ and $MD$ long.
This is impossible in the negative triangle $M$.

Suppose that $w_1w_2$, represented by $OV$, is long. The other edges
of $O,V$ are short, so $R=E$ and $L=F$.
We follow the two arms $R,X,H$ and $L,Y,K$ of the diagram.

\smallskip\noindent\emph{Case 1. $EF$ is long.}
Then $EH$ and $FK$ are short, since a negative triangle has at most
one long edge. The end regions of each arm have opposite orders:
$R\ne H$ and $L\ne K$. Thus $X$ (respectively, $Y$) agrees with
exactly one end. The edge to that end is long, so $XA$ and $YB$
are short.

Consider the upper arm, where $OR$ and $EH$ are already short.
If $HU$ is short, the boundary-edge rules at $H,X,R$ successively
force $HX$ to be long, $RX$ to be short, and $RC$ to be long.
Since $RX$ is short and $R\ne H$, we have $X=H$.
The short edges $XA$ and $HU$ then give $A=U$.
If $HU$ is long, then $HX$ is short, so $X=R$.
Thus $RX$ is long and $RC$ is short; the short edges $XA,RC$ give
$A=C$.
The same argument on the lower arm gives $B=U$ and $LD$ long when
$KU$ is short, and $B=D$ when $KU$ is long.

Since $U$ has exactly one long edge, at least one of $HU,KU$ is short.
Consequently at least one of $RC,LD$ is long, forcing $CD$ to be
short and hence $C\ne D$.
The edges $HU,KU$ cannot both be short, since that would give $A=B=U$,
contrary to Claim~\ref{clm:circular-orders}(ii) at $Q$.
Thus either $A=C$ and $B=U$, or $B=D$ and $A=U$.
The same claim at $Q$ gives $A\ne B$.
Since $C\ne D$, either alternative implies both $A=C$ and $B=D$,
contradicting the central obstruction.

\smallskip\noindent\emph{Case 2. $EF$ is short.}
The exactly-one-long rule at $E,F$ makes $EH,FK$ long.
Consequently $HU,KU$ are short, and the rule at $U$ makes $UQ$ long.
Claim~\ref{clm:circular-orders}(iii) gives $A\ne B$.
We also have $X=Y=U$ and $R=E\ne F=L$, so exactly one of $R,L$
has the same circular order as $U$.
If it is $R$, then $RX$ is long, so $A=C\ne U$. Thus $B=U=Y$,
which makes $YB$ long, $LY$ short and $LD$ long. Therefore $CD$
is short and $B=D$. If it is $L$, the same argument with the two
arms interchanged gives
$B=D\ne U$ and $A=C$. Either way, the central obstruction is violated.
Thus $w_1w_2$ is short.
\end{proof}

\subsection{Proof of the main theorem}\label{subsec:main-proof}

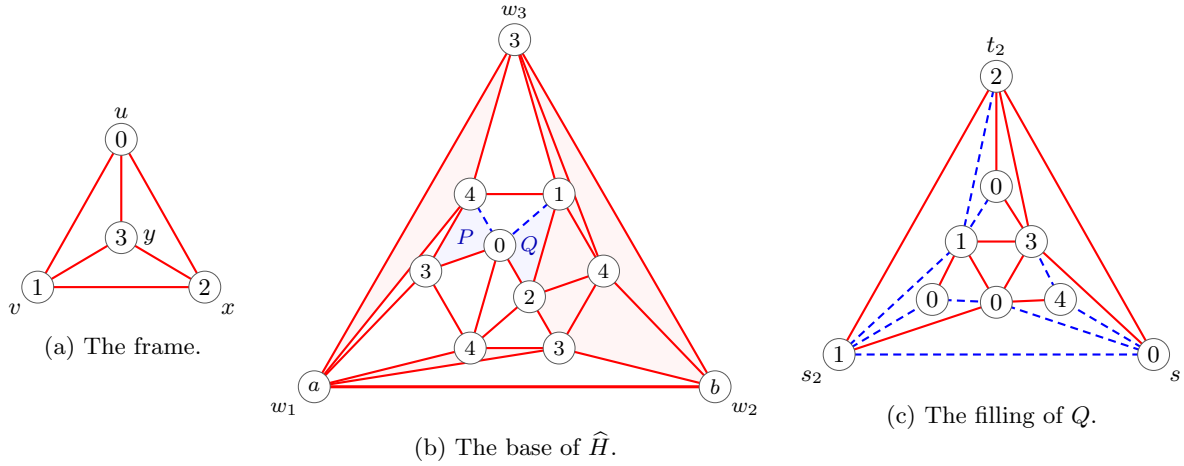
\begin{figure}[htbp]
\centering
\begin{minipage}[c]{.20\textwidth}\centering
\begin{tikzpicture}[x=.85cm,y=.85cm,
  every node/.style={font=\footnotesize,inner sep=1pt},
  negative/.style={draw=red,line width=.8pt},
  vertex/.style={circle,draw=black!65,fill=white,line width=.35pt,minimum size=12pt}]
\coordinate (a) at (0,1.5);
\coordinate (b) at (-1.3,-.8);
\coordinate (c) at (1.3,-.8);
\coordinate (d) at (0,-.03);
\draw[negative] (a)--(b)--(c)--cycle;
\draw[negative] (a)--(d)--(b);
\draw[negative] (c)--(d);
\node[vertex] at (a) {$0$};
\node[vertex] at (b) {$1$};
\node[vertex] at (c) {$2$};
\node[vertex] at (d) {$3$};
\node[above=7pt] at (a) {$u$};
\node[below left=5pt] at (b) {$v$};
\node[below right=5pt] at (c) {$x$};
\node[right=7pt] at (d) {$y$};
\end{tikzpicture}
\par\smallskip\footnotesize (a) The frame.
\end{minipage}\hfill
\begin{minipage}[c]{.44\textwidth}\centering
\begin{tikzpicture}[x=.56cm,y=.56cm,
  every node/.style={font=\scriptsize,inner sep=1pt},
  positive/.style={draw=blue,dash pattern=on 3pt off 2pt,line width=.8pt},
  negative/.style={draw=red,line width=.8pt},
  vertex/.style={circle,draw=black!65,fill=white,line width=.35pt,minimum size=12pt}]
\coordinate (a) at (-4.734272,0.000000);
\coordinate (b) at (0.000000,8.200000);
\coordinate (c) at (4.734272,0.000000);
\coordinate (d) at (-1.050000,4.551987);
\coordinate (e) at (1.050000,4.551987);
\coordinate (f) at (-2.100000,2.733333);
\coordinate (g) at (1.050000,0.914680);
\coordinate (h) at (2.100000,2.733333);
\coordinate (i) at (-1.050000,0.914680);
\coordinate (j) at (-0.350000,3.339551);
\coordinate (k) at (0.350000,2.127116);
\fill[red!4] (a)--(b)--(d)--cycle;
\fill[red!4] (b)--(c)--(h)--cycle;
\fill[red!4] (b)--(e)--(h)--cycle;
\fill[blue!4] (d)--(j)--(f)--cycle;
\fill[blue!4] (e)--(k)--(j)--cycle;
\fill[red!4] (e)--(h)--(k)--cycle;
\fill[red!4] (g)--(h)--(k)--cycle;
\fill[red!4] (c)--(g)--(h)--cycle;
\fill[red!4] (a)--(i)--(g)--cycle;
\draw[negative] (a)--(b);
\draw[negative,line width=1.15pt] (a)--(c);
\draw[negative] (a)--(d);
\draw[negative] (a)--(f);
\draw[negative] (a)--(g);
\draw[negative] (a)--(i);
\draw[negative] (b)--(c);
\draw[negative] (b)--(d);
\draw[negative] (b)--(e);
\draw[negative] (b)--(h);
\draw[negative] (c)--(g);
\draw[negative] (c)--(h);
\draw[negative] (d)--(e);
\draw[negative] (d)--(f);
\draw[positive] (d)--(j);
\draw[negative] (e)--(h);
\draw[positive] (e)--(j);
\draw[negative] (e)--(k);
\draw[negative] (f)--(i);
\draw[negative] (f)--(j);
\draw[negative] (g)--(h);
\draw[negative] (g)--(i);
\draw[negative] (g)--(k);
\draw[negative] (h)--(k);
\draw[negative] (i)--(j);
\draw[negative] (i)--(k);
\draw[negative] (j)--(k);
\node[vertex] at (a) {$a$};
\node[vertex] at (b) {$3$};
\node[vertex] at (c) {$b$};
\node[vertex] at (d) {$4$};
\node[vertex] at (e) {$1$};
\node[vertex] at (f) {$3$};
\node[vertex] at (g) {$3$};
\node[vertex] at (h) {$4$};
\node[vertex] at (i) {$4$};
\node[vertex] at (j) {$0$};
\node[vertex] at (k) {$2$};
\node[above=7pt] at (b) {$w_3$};
\node[below left=5pt] at (a) {$w_1$};
\node[below right=5pt] at (c) {$w_2$};
\node[font=\scriptsize,text=blue!70!black] at (-1.16667,3.54162) {$P$};
\node[font=\scriptsize,text=blue!70!black] at (.35000,3.33955) {$Q$};
\end{tikzpicture}
\par\smallskip\footnotesize (b) The base of $\widehat H$.
\end{minipage}\hfill
\begin{minipage}[c]{.34\textwidth}\centering
\begin{tikzpicture}[x=1.6cm,y=1.6cm,
  every node/.style={font=\footnotesize,inner sep=.8pt},
  positive/.style={draw=blue,dash pattern=on 3pt off 2pt,line width=.8pt},
  negative/.style={draw=red,line width=.8pt},
  vertex/.style={circle,draw=black!65,fill=white,line width=.35pt,minimum size=12pt}]
\coordinate (va) at (-1.300000,-0.800000);
\coordinate (vb) at (0.000000,1.500000);
\coordinate (vc) at (1.300000,-0.800000);
\coordinate (vx) at (0.000000,-0.368196);
\coordinate (vy) at (-0.290000,0.134098);
\coordinate (vz) at (0.290000,0.134098);
\coordinate (vda) at (-0.530000,-0.344699);
\coordinate (vdb) at (0.000000,0.589399);
\coordinate (vdc) at (0.530000,-0.344699);
\draw[negative] (va)--(vb)--(vc);
\draw[positive] (vc)--(va);
\draw[positive] (va)--(vy);
\draw[negative] (vb)--(vz);
\draw[positive] (vc)--(vx);
\draw[negative] (va)--(vx);
\draw[positive] (vb)--(vy);
\draw[negative] (vc)--(vz);
\draw[negative] (vx)--(vy)--(vz)--cycle;
\draw[positive] (va)--(vda)--(vx);
\draw[negative] (vda)--(vy);
\draw[negative] (vb)--(vdb);
\draw[positive] (vdb)--(vy);
\draw[negative] (vdb)--(vz);
\draw[positive] (vc)--(vdc)--(vz);
\draw[negative] (vdc)--(vx);
\node[vertex] at (va) {$1$};
\node[vertex] at (vb) {$2$};
\node[vertex] at (vc) {$0$};
\node[vertex] at (vx) {$0$};
\node[vertex] at (vy) {$1$};
\node[vertex] at (vz) {$3$};
\node[vertex] at (vda) {$0$};
\node[vertex] at (vdb) {$0$};
\node[vertex] at (vdc) {$4$};
\node[above=7pt] at (vb) {$t_2$};
\node[below left=5pt] at (va) {$s_2$};
\node[below right=5pt] at (vc) {$s$};
\end{tikzpicture}
\par\smallskip\footnotesize (c) The filling of $Q$.
\end{minipage}
\caption{A circular $5$-coloring, with colors read modulo $5$.
In (b), take $(a,b)=(0,1),(1,2),(2,0)$ for
$\widehat H_1,\widehat H_2,\widehat H_3$, respectively.
Panel (c) shows $T^+$ after switching its vertex $u$.
To fill $P$, negate the colors in (c), identifying
$(s_2,t_2,s)$ with $(s_1,t,s)$.
Red regions receive $T^-$ fillings as described in the proof.}
\label{fig:five-coloring}
\end{figure}

\begin{proof}[Proof of Theorem~\ref{thm:circular}]
Suppose that $(G,\sigma)$ has a circular $r$-coloring with $r<5$.
Its all-negative $K_4$ requires $r\ge4$.
Lemma~\ref{circ:S} makes $uv,vx,xu$ all short, contradicting
Lemma~\ref{lem:filling}(i) for the all-negative $K_4$.
Hence $\chic(G,\sigma)\ge5$.

For the reverse inequality, use the colors shown in
Figure~\ref{fig:five-coloring}. In each $T^-$ filling, give the added
vertex any color in $\{0,1,2,3,4\}$ absent from its boundary.
Every negative edge then joins distinct integer colors, and every
positive edge joins equal or consecutive colors modulo $5$.
Their circular distances are therefore at least $1$ and at most
$1<5/2-1$, respectively. This is a circular $5$-coloring, proving
$\chic(G,\sigma)=5$.
\end{proof}

\section{Concluding remarks}\label{sec:concluding}

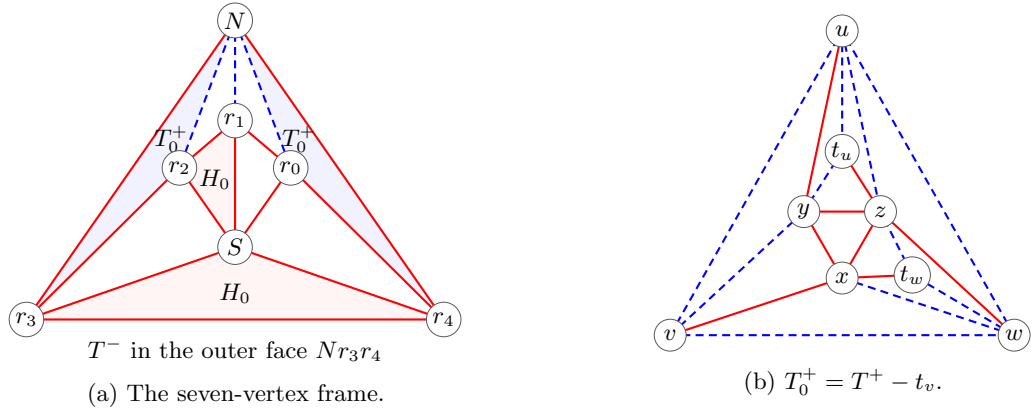
\begin{figure}[H]
\centering
\begin{minipage}[c]{.57\textwidth}
\centering
\begin{tikzpicture}[x=.92cm,y=.92cm,
 every node/.style={font=\footnotesize,inner sep=1pt},
 negative/.style={draw=red,line width=.8pt},
 positive/.style={draw=blue,dash pattern=on 3pt off 2pt,line width=.8pt},
 vertex/.style={circle,draw=black!65,fill=white,line width=.35pt,minimum size=13pt}]
\coordinate (N) at (0,4.3);
\coordinate (r3) at (-3,0);
\coordinate (r4) at (3,0);
\coordinate (S) at (0,1.04);
\coordinate (r2) at (-.80,2.18);
\coordinate (r1) at (0,2.86);
\coordinate (r0) at (.80,2.18);
\fill[blue!5] (N)--(r4)--(r0)--cycle;
\fill[blue!5] (N)--(r2)--(r3)--cycle;
\fill[red!4] (S)--(r3)--(r4)--cycle;
\fill[red!4] (S)--(r1)--(r2)--cycle;
\draw[negative] (r0)--(r1)--(r2)--(r3)--(r4)--cycle;
\foreach \i in {0,1,2}{\draw[positive] (N)--(r\i);}
\foreach \i in {3,4}{\draw[negative] (N)--(r\i);}
\foreach \i in {0,...,4}{\draw[negative] (S)--(r\i);}
\node at (-.92,2.60) {$T^+_0$};
\node at (.92,2.60) {$T^+_0$};
\node at (0,.38) {$H_0$};
\node at (-.29,2.02) {$H_0$};
\foreach \i in {0,...,4}{\node[vertex] at (r\i) {$r_\i$};}
\node[vertex] at (N) {$N$};
\node[vertex] at (S) {$S$};
\node[below] at (0,-.22) {$T^-$ in the outer face $Nr_3r_4$};
\end{tikzpicture}
\par\smallskip\footnotesize (a) The seven-vertex frame.
\end{minipage}\hfill
\begin{minipage}[c]{.41\textwidth}
\centering
\begin{tikzpicture}[x=1.75cm,y=1.75cm,
 every node/.style={font=\footnotesize,inner sep=.8pt},
 negative/.style={draw=red,line width=.8pt},
 positive/.style={draw=blue,dash pattern=on 3pt off 2pt,line width=.8pt},
 vertex/.style={circle,draw=black!65,fill=white,line width=.35pt,minimum size=12pt}]
\coordinate (v0) at (-1.3,-.8);
\coordinate (v1) at (0,1.5);
\coordinate (v2) at (1.3,-.8);
\coordinate (v3) at (0,-.368196);
\coordinate (v4) at (-.29,.134098);
\coordinate (v5) at (.29,.134098);
\coordinate (v7) at (0,.589399);
\coordinate (v8) at (.53,-.344699);
\draw[positive] (v0)--(v1)--(v2)--cycle;
\draw[positive] (v0)--(v4);
\draw[positive] (v1)--(v5);
\draw[positive] (v2)--(v3);
\draw[negative] (v0)--(v3);
\draw[negative] (v1)--(v4);
\draw[negative] (v2)--(v5);
\draw[negative] (v3)--(v4)--(v5)--cycle;
\draw[positive] (v1)--(v7)--(v4);
\draw[negative] (v7)--(v5);
\draw[positive] (v2)--(v8)--(v5);
\draw[negative] (v8)--(v3);
\node[vertex] at (v0) {$v$};
\node[vertex] at (v1) {$u$};
\node[vertex] at (v2) {$w$};
\node[vertex] at (v3) {$x$};
\node[vertex] at (v4) {$y$};
\node[vertex] at (v5) {$z$};
\node[vertex] at (v7) {$t_u$};
\node[vertex] at (v8) {$t_w$};
\end{tikzpicture}
\par\smallskip\footnotesize (b) $T^+_0=T^+-t_v$.
\end{minipage}
\caption{The frame of the $68$-vertex example and its smaller positive
piece. The $H_0$ pieces use the frame of
Figure~\ref{fig:H-construction}, with the two specified $T^-$ fillings
omitted. Red solid edges are negative and blue dashed edges are positive.
Only the four labelled bounded faces and the outer face receive fillings.}
\label{fig:frame68}
\end{figure}

\paragraph{A smaller example.}
A different arrangement gives a $68$-vertex example with circular
chromatic number $5$; its frame is shown in
Figure~\ref{fig:frame68}. Let $H_0$ be obtained from $\widehat H$ by
omitting the $T^-$ fillings in the two triangles $w_3s_2z$ and
$t_3zt_2$, labelled $H$ and $K$ in Figure~\ref{fig:H-construction}.
Thus $H_0$ uses five copies of $T^-$ and two copies of $T^+$ and has
$28$ vertices. Let $T^+_0=T^+-t_v$, with the labels of
Figure~\ref{fig:small-fillings}(b); this graph has eight vertices.

Start with the pentagonal bipyramid with apices $N,S$ and equatorial
cycle $r_0r_1r_2r_3r_4r_0$. All its edges are negative except
$Nr_0,Nr_1,Nr_2$. Insert $H_0$ into the two faces with
\[
 (w_1,w_2,w_3)=(r_4,r_3,S)\quad\text{and}\quad(r_1,r_2,S).
\]
Insert copies of $T^+_0$ with $(u,v,w)=(N,r_0,r_4)$ and $(N,r_2,r_3)$,
switching $w$ in each copy so that its signs match the face boundary.
Fill $Nr_3r_4$ with $T^-$, and leave the other five faces unfilled.
The resulting triangulation has $68$ vertices and $198$ edges.
The lower-bound argument for this smaller graph requires additional
cases, so we have used the $85$-vertex construction for the proof above.

\paragraph{The supremum.}
Theorem~\ref{thm:circular}, together with the upper bound of
Naserasr, Wang, and Zhu~\cite{NWZ2021}, gives
\[
 5\le\chic(\mathcal{SP})\le6.
\]
Determining the supremum $\chic(\mathcal{SP})$ remains open.

\section*{Acknowledgments}
     Lujia Wang is partially supported by the National Natural Science Foundation of China (Grant No. 12371359).

\section*{Contributions and use of AI}
The authors used AI tools from OpenAI to assist in the discovery of the initial 68-vertex example. They subsequently modified the construction to obtain the example presented in this paper, which admits a simpler proof. AI tools were also used to assist with proof checking and manuscript preparation. The authors take full responsibility for all proofs, computations, and conclusions presented in this paper.

\end{document}